\documentclass[leqno,12pt]{article} 
\usepackage{amsmath, amssymb}
\usepackage{amsthm} 
\theoremstyle{plain} 
\newtheorem{theorem}{\indent\sc Theorem}[section]

\newtheorem{corollary}[theorem]{\indent\sc Corollary}
\newtheorem{proposition}[theorem]{\indent\sc Proposition}

\theoremstyle{definition} 
\newtheorem{definition}[theorem]{\indent\sc Definition}
\newtheorem{remark}[theorem]{\indent\sc Remark}

\input{tcilatex}

\title{Bi-$f$-harmonic curves and hypersurfaces} 
\author{Selcen Y\"{u}ksel\ Perkta\c{s}, Adara Monica Blaga, Feyza Esra Erdo\u{g}an,\\
and Bilal Eftal Acet}
\date{}

\begin{document}

\maketitle

\markboth{{\small\it {\hspace{4cm} Bi-$f$-harmonic curves and hypersurfaces}}}{\small\it{Bi-$f$-harmonic curves and hypersurfaces
\hspace{4cm}}}

\footnote{2010 \textit{Mathematics Subject Classification}. 31A30, 53A07, 53C42. }
\footnote{\textit{Key words and phrases}. bi-$f$-harmonic submanifold.}

\begin{abstract} In the present paper, we study bi-$f$-harmonic maps which generalize not only
$f$-harmonic maps, but also biharmonic maps. We derive bi-$f$-harmonic
equations for curves in the Euclidean space, unit sphere, hyperbolic space, and in hypersurfaces of Riemannian manifolds.
\end{abstract}

\section{Introduction}

Harmonic maps between Riemannian manifolds, which can be viewed as a
generalization of geodesics when the domain is $1$-dimensional, or of harmonic
functions when the ranges are Euclidean spaces, have an extensive study area
and there exist many applications of such mappings in mathematics and
physics.
challenge to prove the existence of harmonic maps.
A harmonic map may not
always exist in a homotopy class, and if it exists, then it may not be unique.

As a generalization of harmonic maps, biharmonic maps between Riemannian
manifolds were introduced by J. Eells and J. H. Sampson in \cite%
{Eells-Sampson}.
B. Y. Chen \cite{Chen} defined biharmonic submanifolds of the Euclidean
space
and stated a well-known conjecture: Any biharmonic
submanifold of the Euclidean space is harmonic, thus minimal. If one uses
the definition of biharmonic maps for Riemannian immersions into Euclidean
space, it is easy to see that Chen's definition of biharmonic submanifold
coincides with the definition given by using bienergy functional.
In recent years, there has been achieved an important literature on biharmonic
submanifold theory including many results on the non-existence of biharmonic
submanifolds in manifolds with non-positive sectional curvature. These
non-existence consequences (see \cite{Jiang2}, \cite{MO}) as well as
\textit{Generalized Chen's conjecture}: Any biharmonic submanifold in a Riemannian
manifold with non-positive sectional curvature is minimal, which was
proposed by R. Caddeo, S. Montaldo and C. Oniciuc \cite{CMO1}, led the
studies to spheres and other non-negatively curved spaces. But in recent
years, the authors of \cite{OT} proved that the \textit{Generalized Chen's conjecture} is
not true by constructing examples of proper biharmonic hypersurfaces in a
$5$-dimensional space of non-constant negative sectional curvature. For some
recent geometric studies of general biharmonic maps and biharmonic
submanifolds see (\cite{Ou4}, \cite{CMO1}, \cite{Ou3}, \cite{OT}, \cite{OW}, \cite{SYP1}, \cite{SYP2}, \cite{SYP3}) and the references
therein.

$f$-harmonic maps between Riemannian manifolds were first introduced and
studied by A. Lichnerowicz in 1970 (see also \cite{EL}). They have also some physical meanings by considering
them as solutions of continuous spin systems and inhomogenous Heisenberg
spin systems \cite{p}. Moreover, there is a strong relationship
between $f$-harmonic maps and gradient Ricci solitons \cite{rimoldi}.

There are two ways to formalize such a link between biharmonic maps and $f$%
-harmonic maps. The first formalization is that by mimicking the theory for
biharmonic maps, the authors of \cite{Lu-15} extended bienergy functional to bi-$f$%
-energy functional and obtained a new type of harmonic maps called bi-$f$%
-harmonic maps. This idea was already considered by Ouakkas, Nasri and Djaa
\cite{OND}. They used the terminology \textquotedblleft $f$-biharmonic
maps\textquotedblright\ for the critical points of bi-$f$-energy functional.
As parallel to \textquotedblleft biharmonic maps\textquotedblright , in \cite%
{Lu-15}, they think that it is more reasonable to call them
\textquotedblleft bi-$f$-harmonic maps\textquotedblright. The second
formalization is that by following the definition of $f$-harmonic map, to
extend the $f$-energy functional to the $f$-bienergy functional and obtain another
type of harmonic maps called $f$-biharmonic maps as critical points of $%
f$-bienergy functional.

The concept of $f$-biharmonic maps has been introduced by W.-J. Lu \cite{Lu}
as a generalization of biharmonic maps. A differentiable map between
Riemannian manifolds is said to be $f$-biharmonic if it is a critical point
of the $f$-bienergy functional defined by integral of $f$ times the
square-norm of the tension field, where $f$ is a smooth positive function on
the domain. If $f=1$, then $f$-biharmonic maps are biharmonic. To avoid the
confusion with the types of maps called by the same name in \cite{OND} and
defined as critical points of the square-norm of the $f$-tension field,
some authors (see \cite{Lu}, \cite{Ou5}) called the map defined in \cite{OND}
as \textit{bi-$f$-harmonic map}, which we shall study in this paper.


\section{Preliminaries}

\subsection{Harmonic maps}

\textit{Harmonic maps} $\varphi :(M,g)\rightarrow (N,h)$ between two Riemannian manifolds
are critical points of the energy functional:%
\begin{equation}
E(\varphi )=\frac{1}{2}\int_{\Omega }|d\varphi |^{2}\vartheta
_{g},  \label{sp7}
\end{equation}%
where $\Omega \subset M$ is a compact domain. The corresponding Euler-Lagrange equation
is \cite{Eells-Sampson}:
\begin{equation}
\tau (\varphi )\equiv \func{trace}\nabla d\varphi =0,  \label{sp8}
\end{equation}%
where $\nabla $ is the connection induced from the Levi-Civita
connection $\nabla ^{M}$ of $M$ and the pull-back connection $\nabla
^{\varphi }$. $\tau (\varphi )$ is called the \textit{tension
field} of the map $\varphi$.

\subsection{Biharmonic maps}

\textit{Biharmonic maps} $\varphi :(M,g)\rightarrow (N,h)$ between two Riemannian manifolds
are critical points of the bienergy functional:
\begin{equation}
E_{2}(\varphi )=\frac{1}{2}\int_{\Omega }|\tau (\varphi )|^{2}\vartheta
_{g},
\label{sp9}
\end{equation}%
where $\Omega \subset M$ is a compact domain. The corresponding Euler-Lagrange equation
is \cite{Jiang1}:
\begin{equation}
\tau _{2}(\varphi )\equiv \func{trace}\left(\nabla ^{\varphi }\nabla ^{\varphi }\tau
(\varphi )-\nabla _{\nabla^M
}^{\varphi }\tau
(\varphi )+R^{N}(\tau
(\varphi ),d\varphi)d\varphi \right)=0,  \label{sp10}
\end{equation}%
where $\tau (\varphi )$ is the tension field of $%
\varphi$ and $R^{N}(X,Y):=[\nabla _{X},\nabla _{Y}]-\nabla _{\lbrack X,Y]\text{ }}$ is
the curvature operator on $N$. $\tau_2 (\varphi )$ is called the \textit{bitension
field} of the map $\varphi$.

\smallskip

From the expression of the bitension field $%
\tau _{2}$, it is clear that a harmonic map is automatically a biharmonic
map. So non-harmonic biharmonic maps, which are called proper biharmonic maps,
are more interesting to be studied.

\subsection{$f$-harmonic maps}

\textit{$f$-harmonic maps} $\varphi :(M,g)\rightarrow (N,h)$ between two Riemannian manifolds
are critical points of the $f$-energy functional:%
\begin{equation}
E_{f}(\varphi )=\frac{1}{2}\int_{\Omega }f\,|d\varphi |^{2}\vartheta
_{g},
\label{sp11}
\end{equation}%
where $\Omega \subset M$ is a compact domain. The corresponding Euler-Lagrange equation
is \cite{Co}:
\begin{equation}
\tau _{f}\,(\varphi )\equiv f\,\tau (\varphi )+d\varphi (\func{grad}f)=0,
\label{sp12}
\end{equation}%
where $\tau (\varphi )$ is the tension field of $%
\varphi $. $\tau
_{f}\,(\varphi )$ is called the \textit{$f$-tension field} of the map $\varphi .$

\smallskip

If $f$ is a
constant function, then it is obvious that $f$-harmonic maps are harmonic. So
$f$-harmonic maps, where $f$ is a non-constant function, which are called proper $f$-harmonic maps, are more interesting
to be studied.

\smallskip

There are two ways to formalize a link between biharmonic maps and $f$%
-harmonic maps. Both of them motivate the following definitions.

\subsection{$f$-biharmonic maps}

\textit{$f$-biharmonic maps} $\varphi :(M,g)\rightarrow (N,h)$ between two Riemannian manifolds
are critical points of the $f$-bienergy functional:%
\begin{equation}
E_{2,f}\,(\varphi )=\frac{1}{2}\int_{\Omega }f\,\,|\tau (\varphi )|^{2}\vartheta
_{g},
\label{sp13}
\end{equation}%
where $\Omega \subset M$ is a compact domain. The corresponding Euler-Lagrange equation
is \cite{Lu}:
\begin{equation}
\tau _{2,f}\,(\varphi )\equiv f\,\tau _{2}(\varphi )+(\Delta f)\tau
(\varphi )+2\nabla _{{\func{grad}}f\,}^{\varphi }\,\tau (\varphi )=0,  \label{sp14}
\end{equation}%
where $\tau (\varphi )$ and $\tau _{2}(\varphi )$ are the tension and
bitension fields of $\varphi $, respectively. $\tau _{2,f}\,(\varphi )$
is called the \textit{$f$-bitension field} of the map $\varphi .$

\subsection{Bi-$f$-harmonic maps}

\textit{Bi-$f$-harmonic maps} $\varphi :(M,g)\rightarrow (N,h)$ between two Riemannian manifolds
are critical points of the bi-$f$-energy functional:%
\begin{equation}
E_{f,2}\,(\varphi )=\frac{1}{2}\int_{\Omega }\,|\tau _{f}\,(\varphi
)|^{2}\vartheta
_{g},  \label{sp15}
\end{equation}%
where $\Omega \subset M$ is a compact domain. The corresponding Euler-Lagrange equation
is \cite{OND}:%
\begin{equation}
\tau _{f,2}\,(\varphi )\equiv -\func{trace} \left(\nabla ^{\varphi
}f\,\left( \nabla ^{\varphi }\tau _{f}\,\left( \varphi \right) \right)
-f\,\nabla _{\nabla ^{M}}^{\varphi }\tau _{f}\,\left( \varphi \right)+f\,R^{N}\left(\tau _{f\,}\left( \varphi \right),  d\varphi \right)
d\varphi\right)=0,  \label{sp16}
\end{equation}%
where $\tau_f (\varphi )$ is the $f$-tension field of $%
\varphi$. $\tau _{f,2}\,(\varphi )$
is called the \textit{bi-$f$-tension field} of the map $\varphi .$

\smallskip

The following inclusions illustrate the relations among these different types of
harmonic maps:%
\[
\text{harmonic maps}\subset \text{biharmonic maps}\subset f\text{%
-biharmonic maps,}
\]%
\[
\text{harmonic maps}\subset f\text{-harmonic maps}\subset \text{%
bi-}f\text{-harmonic maps.}
\]

\section{{\protect\Large Bi-}${\protect\Large f}${\protect\Large -harmonic
curves }}

In this section we derive the bi-$f$-harmonic equation for curves in Riemannian
manifolds and discuss the particular case of the Euclidean space, unit sphere and hyperbolic space.
The following proposition for Euler-Lagrange
equation of bi-$f$-harmonic maps originates from \cite{OND}.

\begin{proposition}
Let $\varphi :(M,g)\rightarrow (N,h)$ be a smooth map between Riemannian
manifolds. Then, in terms of Euler-Lagrange equation, $\varphi $ is a bi-$f$%
-harmonic map if and only if its bi-$f$-tension field $\tau _{f,2}(\varphi )$ vanishes, i.e.
\begin{equation}
\func{trace} \left(\nabla ^{\varphi
}f\,\left( \nabla ^{\varphi }\tau _{f}\,\left( \varphi \right) \right)
-f\,\nabla _{\nabla ^{M}}^{\varphi }\tau _{f}\,\left( \varphi \right)+f\,R^{N}\left( \tau _{f\,}\left( \varphi \right), d\varphi \right)
d\varphi\right)=0,  \label{f-bi-1}
\end{equation}%
where $f :I\rightarrow (0,\infty)$ is a smooth map defined on a real interval $I$ and $\tau _{f}\,\left( \varphi \right) $ is the $f$-tension field given by (%
\ref{sp12}).
\end{proposition}

Clearly, it is observed from (\ref{f-bi-1}) that bi-$f$-harmonic map is a
much wider generalization of harmonic map, because it is not only a
generalization of $f$-harmonic map (as $f\neq 1$ and $\tau _{f}(\varphi )=0$%
), but also a generalization of biharmonic map (as $f=1$). Therefore, it
would be interesting to know whether there is any non-trivial or proper bi-$%
f$-harmonic map which is neither harmonic map nor $f$-harmonic map with $%
f\neq $ constant.

\begin{definition}
A submanifold in a Riemannian manifold is called a \textit{bi-$f$-har\-mo\-nic
submanifold} if the isometric immersion defining the submanifold is a bi-$f$%
-har\-mo\-nic map.
\end{definition}

\bigskip

Let $\gamma :I\rightarrow (N,h)$ be
a curve in a Riemannian manifold $(N,h)$, defined on an open
real interval $I$ and parametrized by its arclength, and $%
\gamma ^{\prime }=:T.$ We have%
\[
\tau \,\left( \gamma \right) =\nabla _{T}^{N}T
\]%
\[
\tau _{f}\,\left( \gamma \right) =f\,\nabla _{T}^{N}T+f^{\,\prime }T
\]%
and in order to obtain the bi-$f$-tension field of $\gamma $, we compute:
\begin{eqnarray}
\func{trace}\left( \nabla ^{\gamma }f\,\left( \nabla ^{\gamma }\tau _{f\,\,}\left(
\gamma \right) \right) -f\,\nabla _{\nabla ^{M}}^{\gamma }\tau _{f}\,\left(
\gamma \right) \right) &=&\nabla _{\frac{d}{dt}}^{\gamma }f\,\left( \nabla
_{\frac{d}{dt}}^{\gamma }\tau _{f}\left( \gamma \right) \right) -f\,\nabla
_{\nabla _{\frac{d}{dt}}^{I}\frac{d}{dt}}^{\gamma }\tau _{f}\left( \gamma
\right)  \nonumber \\
&=&\nabla _{T}^{N}f\,( \nabla _{T}^{N}(f\,\nabla _{T}^{N}T+f^{\,\prime
}T) )  \nonumber \\
&=&\left( ff^{\,\prime \prime \prime
}+f^{\,\prime }f^{\,\prime \prime }\right) T+( 3ff^{\,\prime \prime
}+2\left( f^{\,\prime }\right) ^{2}) \nabla _{T}^{N}T  \nonumber \\
&&+4ff^{\,\prime }\,\nabla _{T}^{N}\nabla _{T}^{N}T+f^{2}\nabla
_{T}^{N}\nabla _{T}^{N}\nabla _{T}^{N}T  \label{3}
\end{eqnarray}%
and
\begin{eqnarray}
\func{trace}\left( R^{N}\left( \tau _{f}\left( \gamma \right), d\gamma \right)
d\gamma \right)&=&R^{N}\left( \tau _{f}\left( \gamma \right), d\gamma\left( \frac{d}{dt}\right)
 \right) d\gamma \left( \frac{d}{dt}\right)
\nonumber \\
&=&fR^{N}\left( \nabla _{T}^{N}T, T\right) T.  \label{2a}
\end{eqnarray}%

From (\ref{3}) and (\ref{2a}) we obtain


\begin{proposition}
~\noindent Let $\gamma :I\rightarrow (N,h)$ be a
curve in a Riemannian manifold $(N,h)$, parametrized by its arclength, and $\gamma
^{\prime }=T.$ Then $\gamma $ is a bi-$f$-harmonic curve if and only if
\begin{eqnarray}
0 &=&\left( ff^{\,\prime \prime \prime
}+f^{\,\prime }f^{\,\prime \prime }\right) T+( 3ff^{\,\prime \prime
}+2\left( f^{\,\prime }\right) ^{2}) \nabla _{T}^{N}T  \nonumber \\
&&+4ff^{\,\prime }\,\nabla _{T}^{2}T+f^{\,2}\nabla
_{T}^{3}T+f^{\,2}R^{N}\left( \nabla _{T}^{N}T, T\right) T,  \label{f-bi-3}
\end{eqnarray}%
where $f:I\rightarrow (0,\infty)$ is a smooth map, $\nabla _{T}^{2}T=:\nabla _{T}^{N}\nabla
_{T}^{N}T$ and $\nabla _{T}^{3}T=:\nabla _{T}^{N}\nabla _{T}^{N}\nabla
_{T}^{N}T.$
\end{proposition}

\noindent \qquad

Let $\left\{ E_{1},E_{2},...,E_{n}\right\}$ be the Frenet frame on the $n$-dimensional manifold $N$,
defined along $\gamma,$ where $E_{1}=\gamma ^{\prime }=T$ is the unit
tangent vector field of $\gamma $, $E_{2}$ is the unit normal vector field
of $\gamma,$ with the same direction as $\nabla _{T}^{N}E_{1}$ and the vector fields $%
E_{3},...,E_{n}$ are the unit vector fields obtained from the Frenet equations for
$\gamma:$%
\begin{equation}
\left\{
\begin{array}{c}
\nabla_{T}^{N}E_{1} =k_{1}E_{2}, \\
\nabla_{T}^{N}E_{2} =-k_{1}E_{1}+k_{2}E_{3}, \\
... \\
\nabla_{T}^{N}E_{r} =-k_{r-1}E_{r-1}+k_{r}E_{r+1},\quad r=3,...,n-1, \\
... \\
\nabla_{T}^{N}E_{n} =-k_{n-1}E_{n-1},
\end{array}%
\right.   \label{f-bi-4}
\end{equation}%
where $k_{1}=\left\Vert \nabla _{T}^{N}E_{1}\right\Vert $ and $k_{2},...,k_{n-1}$
are real valued non-negative maps.


From (\ref{f-bi-4}) we have

\begin{equation}
\nabla _{T}^{2}T=\nabla _{T}^{N}\nabla
_{T}^{N}T=-k_{1}^{2}E_{1}+k_{1}^{\prime }E_{2}+k_{1}k_{2}E_{3},  \label{5}
\end{equation}%
\begin{eqnarray}
\nabla _{T}^{3}T &=&\nabla _{T}^{N}\nabla _{T}^{N}\nabla _{T}^{N}T  \nonumber
\\
&=&-3k_{1}k_{1}^{\prime }E_{1}+\left( k_{1}^{\prime \prime
}-k_{1}^{3}-k_{1}k_{2}^{2}\right) E_{2} \nonumber  \\
&&+\left( 2k_{1}^{\prime }k_{2}+k_{1}k_{2}^{\prime }\right)
E_{3}+k_{1}k_{2}k_{3}E_{4},  \label{6a}
\end{eqnarray}%
\begin{equation}
R^{N}(\nabla _{T}^{N}T, T)T=k_{1}R^{N}(E_2,E_{1})E_1.  \label{7}
\end{equation}%

Using (\ref{f-bi-4}), (\ref{5}), (\ref{6a}) and (\ref{7}) in (\ref{f-bi-3}), we have




\begin{theorem}
\noindent Let $\gamma :I\rightarrow (N,h)$ be a curve in a Riemannian manifold $(N,h)$, parametrized by
its arclength. Then $\gamma $ is a bi-$f$-harmonic curve if and only if%
\begin{eqnarray}
0 &=&\left( -3k_{1}k_{1}^{\prime }f^{2}-4k_{1}^{2}f\,f^{\,\prime
}+f\,f^{\,\prime \prime \prime }+f^{\,\prime }f^{\,\prime \prime }\right) E_1
\nonumber \\
&&+\left( -k_{1}^{3}f^{2}-k_{1}k_{2}^{2}f^{2}+k_{1}^{\prime
\prime }f^{2}+4k_{1}^{\prime }ff^{\,\prime }+3k_{1}ff^{\,\prime \prime
}+2k_{1}(f^{\,\prime })^{2}\right) E_{2}  \nonumber \\
&&+\left(\left( 2k_{1}^{\prime }k_{2}f+k_{1}k_{2}^{\prime
}f+4k_{1}k_{2}f^{\,\prime }\right) f\right) E_{3} \nonumber  \\
&&+\left( k_{1}k_{2}k_{3}f^{\,2}\right) E_{4}+k_{1}f^{2}R^{N}(E_2,E_{1})E_1.
\label{7aa}
\end{eqnarray}
\end{theorem}

\begin{remark}
The property of a curve of being bi-$f$-harmonic in an $n$-dimensional space (with $n>3$) does not depend on all its curvatures, but only on $k_1$, $k_2$ and $k_3$.
\end{remark}

\bigskip

It is well known that in a Riemannian manifold $(N,h)$ of constant sectional
curvature $c$, the curvature tensor field $R^{N}$ is of the form%
\[
R^{N}(X,Y)Z=c\left( h(Y,Z)X-h(X,Z)Y\right),
\]%
for any $X,Y,Z\in \Gamma (TN).$

\bigskip

Then we have

\begin{theorem}
\noindent Let $\gamma :I\rightarrow (N(c),h)$ be a curve in a Riemannian space form $(N(c),h)$, parametrized by
its arclength. Then $\gamma $ is a bi-$f$-harmonic curve if and only if%
\begin{equation}
\left\{
\begin{array}{c}
-3k_{1}k_{1}^{\prime }f^{\,2}-4k_{1}^{2}f\,f^{\,\prime }+f\,f^{\,\prime
\prime \prime }+f^{\,\prime }f^{\,\prime \prime }=0, \\
-k_{1}^{3}f^{\,2}-k_{1}k_{2}^{2}f^{\,2}+k_{1}^{\prime \prime
}f^{\,2}+4k_{1}^{\prime }ff^{\,\prime }+3k_{1}ff^{\,\prime \prime
}+2k_{1}(f^{\,\prime })^{2}+ck_{1}f^{\,2}=0, \\
2k_{1}^{\prime }k_{2}f+k_{1}k_{2}^{\prime
}f+4k_{1}k_{2}f^{\,\prime }=0, \\
k_{1}k_{2}k_{3}=0.%
\end{array}%
\right.  \label{8a}
\end{equation}
\end{theorem}

\bigskip


Let $\gamma :I\rightarrow \mathbb{E}^n$ be a curve in the $n$-dimensional Euclidean
space, defined on an open real interval $I$ and parametrized by its arclength. Since $\mathbb{E}^n$ is a Riemannian space form with $c=0$, from the bi-$f$%
-harmonic curve equation given by (\ref{8a}) we have

\begin{theorem}
Let $\gamma :I\rightarrow \mathbb{E}^{n}$ be a curve in the
$n$-dimensional Euclidean space, parametrized by its arclength. Then $\gamma $ is a bi-$f$-harmonic
curve if and only if%
\begin{equation}
\left\{
\begin{array}{c}
-3k_{1}k_{1}^{\prime }f^2-4k_{1}^{2}f\,f^{\,\prime }+f\,f^{\,\prime
\prime \prime }+f^{\,\prime }f^{\,\prime \prime }=0, \\
-k_{1}^{3}f^{\,2}-k_{1}k_{2}^{2}f^{\,2}+k_{1}^{\prime \prime
}f^{\,2}+4k_{1}^{\prime }ff^{\,\prime }+3k_{1}ff^{\,\prime \prime
}+2k_{1}(f^{\,\prime })^{2}=0, \\
2k_{1}^{\prime }k_{2}f+k_{1}k_{2}^{\prime
}f+4k_{1}k_{2}f^{\,\prime } =0, \\
k_{1}k_{2}k_{3}=0.%
\end{array}%
\right.  \label{9a}
\end{equation}
\end{theorem}

\bigskip

\qquad \textbf{CASE I: }If $k_{1}=0$, namely $\gamma $ is a geodesic curve,
then from (\ref{9a}) we obtain that it is bi-$f$-harmonic if and only if $ff^{\prime \prime
}=constant.$

\smallskip

It is well known that geodesics are $f$-harmonic with $%
f=constant$ and they are automatically bi-$f$-harmonic$.$ Remark that also for $f(t)=at+b,$ $%
a,b\in \mathbb{R}$, any geodesic curve is bi-$f$-harmonic.



\begin{theorem}
A geodesic curve is bi-$f$-harmonic if and only if $ff^{\prime \prime
}=constant.$
\end{theorem}


\bigskip

\qquad \textbf{CASE II:} \textbf{\ }If $k_{1}=constant\neq 0$ and $%
k_{2}=0, $ then (\ref{9a}) reduces to%
\begin{equation}
\left\{
\begin{array}{c}
-4k_{1}^{2}ff^{\prime }+ff^{\prime
\prime \prime }+f^{\prime }f^{\prime \prime } =0, \\
-k_{1}^{2}f^{2}+3ff^{\prime \prime } +2(f^{\prime })^{2}=0.%
\end{array}%
\right.  \label{9b}
\end{equation}%
From the second equation above we obtain
\begin{equation}
f f^{\prime \prime }=\frac{k_{1}^{2}f^{2}-2(f^{\prime })^{2}}{3},
\label{9c}
\end{equation}%
which implies
\begin{equation}
f^{\prime }\left( 5k_{1}^{2}f+2f^{\prime \prime }\right) =0,  \label{9d}
\end{equation}%
via the first equation of (\ref{9b}) and we get

\begin{theorem}
Let $\gamma :I\rightarrow \mathbb{E}^{n}$ be a curve in the
$n$-dimensional Euclidean space, parametrized by its arclength, with $k_{1}=$ constant $\neq 0$ and $k_{2}=0.$
Then $\gamma $ is a bi-$f$-harmonic curve if and only if either $f$ is a
constant function or $f$ is given by%
\[
f(s)=c_{1}\cos \left( \sqrt{\frac{5}{2}}k_{1}s\right) +c_{2}\sin \left(
\sqrt{\frac{5}{2}}k_{1}s\right) ,
\]%
for $s\in I$ and $c_1,c_2\in \mathbb{R}$.
\end{theorem}

\bigskip

\qquad \textbf{CASE III:} \textbf{\ }If $k_{1}=constant\neq 0$ and $%
k_{2}=constant\neq 0,$ then (\ref{9a}) reduces to%
\begin{equation}
\left\{
\begin{array}{c}
-4k_{1}^{2}ff^{\prime }+ff^{\prime \prime \prime }+f^{\prime }f^{\prime
\prime }=0, \\
-k_{1}^{2}f^{2}-k_{2}^{2}f^{2}+3ff^{\prime \prime }+2(f^{\prime })^{2}=0, \\
f^{\prime }=0, \\
k_{3}=0,%
\end{array}%
\right.  \label{9e}
\end{equation}%
which implies
\begin{equation}
\left\{
\begin{array}{c}
k_{1}^{2}+k_{2}^{2} =0, \\
f^{\prime }=0, \\
k_{3}=0,%
\end{array}%
\right.  \label{9f}
\end{equation}%
and we deduce

\begin{theorem}
There is no bi-$f$-harmonic curve in the
$n$-dimensional Euclidean space with $k_{1}=$ constant $\neq 0$ and $%
k_{2}=constant\neq 0.$
\end{theorem}

\bigskip

\qquad \textbf{CASE IV:} \textbf{\ }If $k_{1}=constant\neq 0$ and $k_{2}\neq constant$, then (\ref{9a}) reduces to%
\begin{equation}
\left\{
\begin{array}{c}
-4k_{1}^{2}ff^{\prime }+ff^{\prime \prime \prime
}+f^{\prime }f^{\prime \prime }=0, \\
-k_{1}^{2}f^{2}-k_{2}^{2}f^{2}+3ff^{\prime \prime }+2(f^{\prime })^{2}=0, \\
k_{2}^{\prime }f+4k_{2}f^{\prime } =0, \\
k_2k_{3}=0,%
\end{array}%
\right.  \label{...}
\end{equation}%
and we have

\begin{theorem}
Let $\gamma :I\rightarrow \mathbb{E}^{n}$ be a curve in
the $n$-dimensional Euclidean space, parametrized by its arclength, with $k_{1}=constant\neq 0$ and $k_{2}\neq constant$ and nowhere zero.
Then $\gamma $ is a bi-$f$-harmonic curve if and
only if $f=ck_{2}^{-\frac{1}{4}}$ (with $c$ a positive constant), $k_{3}=0$ and
the curvatures $k_{1}$ and $k_2$ satisfy:%
\begin{equation}
\left\{
\begin{array}{c}
32 k_1^2k_2^2 k_{2}^{\prime }-25 (k_{2}^{\prime })^3+32 k_2k_{2}^{\prime }k_{2}^{\prime \prime}-8 k_2^2k_{2}^{\prime \prime \prime}=0, \\
16 k_1^2k_2^2+16 k_2^4-17 (k_{2}^{\prime })^2+12 k_2k_{2}^{\prime \prime}=0.%
\end{array}%
\right.  \label{...}
\end{equation}
\end{theorem}

\bigskip

\qquad \textbf{CASE V:} \textbf{\ }Concerning the case $k_{1}\neq constant$ and $k_{2}=0$, we can state

\begin{theorem}
Let $\gamma :I\rightarrow \mathbb{E}^{n}$ be a curve in
the $n$-dimensional Euclidean space, parametrized by its arclength, with $k_{1}\neq constant$ and $k_{2}=0$.
Then $\gamma $ is a bi-$f$-harmonic curve if and
only if the curvatures $k_{1}$ and $k_2$ satisfy:%
\begin{equation}
\left\{
\begin{array}{c}
-3k_{1}k_{1}^{\prime }f^2-4k_{1}^{2}f\,f^{\,\prime }+f\,f^{\,\prime
\prime \prime }+f^{\,\prime }f^{\,\prime \prime }=0, \\
-k_{1}^{3}f^{\,2}+k_{1}^{\prime \prime
}f^{\,2}+4k_{1}^{\prime }ff^{\,\prime }+3k_{1}ff^{\,\prime \prime
}+2k_{1}(f^{\,\prime })^{2}=0.%
\end{array}%
\right.  \label{...}
\end{equation}%
\end{theorem}

\bigskip

\qquad \textbf{CASE VI:} \textbf{\ }If $k_{1}\neq constant$ and $%
k_{2}=constant\neq 0,$ then (\ref{9a}) reduces to%
\begin{equation}
\left\{
\begin{array}{c}
-3k_{1}k_{1}^{\prime }f^{2}-4k_{1}^{2}ff^{\prime }+ff^{\prime \prime \prime
}+f^{\prime }f^{\prime \prime }=0, \\
-k_{1}^{3}f^{2}-k_{1}k_{2}^{2}f^{2}+k_{1}^{\prime \prime
}f^{2}+4k_{1}^{\prime }ff^{\prime }+3k_{1}ff^{\prime \prime }+2k_{1}(f^{\prime })^{2}=0, \\
k_{1}^{\prime }f+2k_{1}f^{\prime } =0, \\
k_1k_{3}=0,%
\end{array}%
\right.  \label{9g}
\end{equation}%
and we have

\begin{theorem}
Let $\gamma :I\rightarrow \mathbb{E}^{n}$ be a curve in
the $n$-dimensional Euclidean space, parametrized by its arclength, with $k_{1}\neq constant$ and nowhere zero and $%
k_{2} = constant\neq 0.$ Then $\gamma $ is a bi-$f$-harmonic curve if and
only if $f=ck_{1}^{-\frac{1}{2}}$ (with $c$ a positive constant), $k_{3}=0$ and
the curvatures $k_{1}$ and $k_2$ satisfy:%
\begin{equation}
\left\{
\begin{array}{c}
9(k_{1}^{\prime })^{3}+4k_{1}^{4}k_{1}^{\prime }-10k_{1}k_{1}^{\prime }k_{1}^{\prime \prime
}+2k_{1}^{2}k_{1}^{\prime \prime \prime }=0, \\
3(k_{1}^{\prime })^{2}-4k_{1}^{4}-4k_{1}^{2}k_{2}^{2}-2k_{1}k_{1}^{\prime \prime
}=0.%
\end{array}%
\right.  \label{9i}
\end{equation}
\end{theorem}

\bigskip

\qquad \textbf{CASE VII:} \textbf{\ }Concerning the case $k_{1}\neq constant$ and $k_{2}\neq
constant,$ we can state

\begin{theorem}
Let $\gamma :I\rightarrow \mathbb{E}^{n}$ be a curve in the
$n$-dimensional Euclidean space, parametrized by its arclength, with $k_{1}\neq constant$ and $%
k_{2}\neq constant$ and $k_1$, $k_2$ are nowhere zero. Then $\gamma $ is a bi-$f$-harmonic curve if and only
if $f=ck_{1}^{-\frac{1}{2}}k_{2}^{-\frac{1}{4}}$ (with $c$ a positive constant), $%
k_{3}=0$ and the curvatures $k_{1}$ and $k_{2}$ satisfy:%
\begin{equation}
\left\{
\begin{array}{c}
-3k_{1}k_{1}^{\prime }f^{2}-4k_{1}^{2}ff^{\prime }+ff^{\prime \prime \prime
}+f^{\prime }f^{\prime \prime }=0, \\
-k_{1}^{3}f^{2}-k_{1}k_{2}^{2}f^{2}+k_{1}^{\prime \prime }f^{2}
+4k_{1}^{\prime }ff^{\prime }+3k_{1}ff^{\prime \prime }+2k_{1}(f^{\prime
})^{2}=0.%
\end{array}%
\right.  \label{10b}
\end{equation}
\end{theorem}

\bigskip

Similar results hold for bi-$f$-harmonic curves in the $n$-dimensional sphere $S^{n}(1)$ and in the $n$-dimensional hyperbolic space $H^{n}(-1)$.

\begin{theorem}
Let $\gamma :I\rightarrow S^{n}(1)$ be a curve parametrized by its arclength. Then $%
\gamma $ is a bi-$f$-harmonic curve if and only if%
\begin{equation}
\left\{
\begin{array}{c}
-3k_{1}k_{1}^{\prime }f^{2}-4k_{1}^{2}ff^{\prime }+ff^{\prime \prime \prime
}+f^{\prime }f^{\prime \prime }=0, \\
-k_{1}^{3}f^{2}-k_{1}k_{2}^{2}f^{2}+k_{1}^{\prime \prime }f^{2}+k_{1}f^{2}
+4k_{1}^{\prime }ff^{\prime }+3k_{1}ff^{\prime \prime }+2k_{1}(f^{\prime
})^{2}=0, \\
2k_{1}^{\prime }k_{2}f+k_{1}k_{2}^{\prime }f+4k_{1}k_{2}f^{\prime
}=0, \\
k_1k_2k_3=0.
\end{array}%
\right.  \label{11a}
\end{equation}
\end{theorem}

\begin{theorem}
Let $\gamma :I\rightarrow H^{n}(-1)$ be a curve parametrized by its arclength. Then $%
\gamma $ is a bi-$f$-harmonic curve if and only if%
\begin{equation}
\left\{
\begin{array}{c}
-3k_{1}k_{1}^{\prime }f^{2}-4k_{1}^{2}ff^{\prime }+ff^{\prime \prime \prime
}+f^{\prime }f^{\prime \prime }=0, \\
-k_{1}^{3}f^{2}-k_{1}k_{2}^{2}f^{2}+k_{1}^{\prime \prime }f^{2}-k_{1}f^{2}
+4k_{1}^{\prime }ff^{\prime }+3k_{1}ff^{\prime \prime }+2k_{1}(f^{\prime
})^{2}=0, \\
2k_{1}^{\prime }k_{2}f+k_{1}k_{2}^{\prime }f+4k_{1}k_{2}f^{\prime
}=0, \\
k_1k_2k_3=0.
\end{array}%
\right.  \label{11a}
\end{equation}
\end{theorem}

Concerning the CASES IV--VII, we obtain similar conditions like in the Euclidean space and in the CASES I--III, we get
the following characterizations of bi-$f$-harmonic curves in $S^n(1)$ and $H^n(-1)$, respectively.

\begin{theorem}
Let $\gamma :I\rightarrow N$ be a curve in
$N$, parametrized by its arclength.
\begin{enumerate}
  \item For $N:=S^n(1)$:
  \begin{enumerate}
    \item if $k_1=0$, then $\gamma $ is a bi-$f$-harmonic curve if and
only if $ff^{\prime \prime
}=constant;$
    \item if $k_{1}=$ constant $\neq 0$ and $k_{2}=0,$
then $\gamma $ is a bi-$f$-harmonic curve if and only if either $f$ is a
constant function or $f$ is given by%
\[
f(s)=c_{1}\cos \left( \sqrt{\frac{5k_{1}^2+1}{2}}s\right) +c_{2}\sin \left(
\sqrt{\frac{5k_{1}^2+1}{2}}s\right) ,
\]%
for $s\in I$ and $c_1,c_2\in \mathbb{R}$;
    \item if $k_{1}=$ constant $\neq 0$ and $%
k_{2}=constant\neq 0,$ then $\gamma $ is a bi-$f$-harmonic curve if and only if $f$ is a constant function, $k_1^2+k_2^2=1$ and $k_3=0$.
  \end{enumerate}
  \item For $N:=H^n(-1)$:
  \begin{enumerate}
    \item if $k_1=0$, then $\gamma $ is a bi-$f$-harmonic curve if and
only if $ff^{\prime \prime
}=constant;$
    \item if $k_{1}=$ constant $\neq 0$ and $k_{2}=0,$
then $\gamma $ is a bi-$f$-harmonic curve if and only if either $f$ is a
constant function or $f$ is given by one of the following expressions%
\[
f(s)=c_{1}s+c_{2}, \ \ \textit{for} \ \  k_1=\pm\frac{\sqrt{5}}{5},
\]%
or
\[
f(s)=c_{1}\cos \left( \sqrt{\frac{5k_{1}^2-1}{2}}s\right) +c_{2}\sin \left(
\sqrt{\frac{5k_{1}^2-1}{2}}s\right), \]%
\[\textit{for} \ \  k_1\in \left(-\infty, -\frac{\sqrt{5}}{5}\right)\cup \left(\frac{\sqrt{5}}{5}, \infty\right),
\]%
or
\[
f(s)=c_{1}e^{\sqrt{\frac{1-5k_{1}^2}{2}}s} +c_{2}e^{
-\sqrt{\frac{1-5k_{1}^2}{2}}s}, \ \ \textit{for} \ \  k_1\in \left(-\frac{\sqrt{5}}{5},\frac{\sqrt{5}}{5}\right),
\]%
for $s\in I$ and $c_1,c_2\in \mathbb{R}$;
    \item if $k_{1}=$ constant $\neq 0$ and $%
k_{2}=constant\neq 0,$ then there is no bi-$f$-harmonic curve.
\end{enumerate}
\end{enumerate}
\end{theorem}

\section{Bi-$f$-harmonic hypersurfaces}

In this section we derive the bi-$f$-harmonic equation for hypersurfaces in
Riemannian manifolds. Let $M$ be an $m$-dimensional hypersurface of $(N,h)$
with mean curvature vector $\eta =H\xi $, where $\xi $ is the unit normal
vector field of $M.$
Denoting also by $h$ the Riemannian metric induced on $M$, by $\nabla^M$ and $\nabla^N$ the Levi-Civita connections on $(M,h)$ and $(N,h)$ respectively, the Gauss and Weingarten formulas corresponding to $M$ are given by:
\begin{equation}
\nabla_{X}^NY=\nabla^M_XY+B(X,Y),
\end{equation}
\begin{equation}
\nabla_{X}^N \xi=-AX,
\end{equation}
for any $X$, $Y\in \Gamma (TM)$, where $B$ is the (symmetric) second fundamental tensor corresponding to $\xi$,
$A$ is the shape operator with respect to the unit normal vector field $\xi$,
and let $b(X,Y)=\langle B(X,Y),\xi \rangle$, for any $X$, $Y\in \Gamma (TM)$.

The bi-$f$-tension field of the immersion $\varphi :M\rightarrow
N$ is given by \cite{OND}:%
\begin{equation}
\tau _{f,2}(\varphi )=-f\func{trace}(\nabla ^{\varphi })^{2}\tau _{f\,}\left(
\varphi \right) -f\func{trace}R^{N}\left( \tau _{f\,}\left( \varphi \right)
,d\varphi \right) d\varphi -\nabla _{\func{grad}f}^{\varphi }\tau
_{f}\,\left( \varphi \right) .  \label{hp2}
\end{equation}%

For an orthonormal frame field $\left\{
e_{1},e_{2},...,e_{m}\right\} \subset \Gamma (TM)$, we have
\begin{eqnarray}
\func{trace}(\nabla ^{\varphi }f\,\left( \nabla ^{\varphi }\tau _{f}\,\left(
\varphi \right) \right) -f\,\nabla _{\nabla ^{M}}^{\varphi }\tau
_{f}\,\left( \varphi \right) ) &=&f\sum\limits_{i=1}^{m}\left\{ \nabla
_{e_{i}}^{\varphi }\nabla _{e_{i}}^{\varphi }\tau _{f\,}\left( \varphi
\right) -\nabla _{\nabla _{e_{i}}^{M}e_{i}}^{\varphi }\tau _{f}\,\left(
\varphi \right) \right\}  \nonumber \\
&&+\nabla _{\func{grad}f}^{\varphi }\tau _{f}\,\left( \varphi \right). \label{hp1}
\end{eqnarray}%

As a first step, we compute
\begin{eqnarray}
\func{trace}(\nabla ^{\varphi })^{2}\tau _{f\,}\left( \varphi \right)
&=&\sum\limits_{i=1}^{m}\left\{ \nabla _{e_{i}}^{\varphi }\nabla
_{e_{i}}^{\varphi }\tau _{f\,}\left( \varphi \right) -\nabla _{\nabla
_{e_{i}}^{M}e_{i}}^{\varphi }\tau _{f}\,\left( \varphi \right) \right\}
\nonumber \\
&=&\sum\limits_{i=1}^{m}\left\{ \nabla _{e_{i}}^{\varphi }\nabla
_{e_{i}}^{\varphi }\left( f\,\tau (\varphi )+d\varphi (\func{grad}f)\right) -\nabla
_{\nabla _{e_{i}}^{M}e_{i}}^{\varphi }\left( f\,\tau (\varphi )+d\varphi
(\func{grad}f)\right) \right\}  \nonumber \\
&=&\sum\limits_{i=1}^{m}\left\{
\begin{array}{c}
\nabla _{e_{i}}^{\varphi }\nabla _{e_{i}}^{\varphi }\left( f\,\tau (\varphi
)\right)+\nabla _{e_{i}}^{\varphi }\nabla _{e_{i}}^{\varphi }d\varphi (\func{grad}f)\\
-\nabla _{\nabla _{e_{i}}^{M}e_{i}}^{\varphi }(f\,\tau (\varphi ))-\nabla
_{\nabla _{e_{i}}^{M}e_{i}}^{\varphi }d\varphi (\func{grad}f)%
\end{array}%
\right\}  \nonumber \\
&=&\sum\limits_{i=1}^{m}\left\{
\begin{array}{c}
\nabla _{e_{i}}^{\varphi }\left( e_{i}(f)\tau (\varphi )+f\nabla
_{e_{i}}^{\varphi }\,\tau (\varphi )\right) +\nabla _{e_{i}}^{\varphi
}\nabla _{e_{i}}^{\varphi }d\varphi (\func{grad}f) \\
-\left( \nabla _{e_{i}}^{M}e_{i}\right) (f)\tau (\varphi )-f\nabla _{\nabla
_{e_{i}}^{M}e_{i}}^{\varphi }\tau (\varphi )-\nabla _{\nabla
_{e_{i}}^{M}e_{i}}^{\varphi }d\varphi (\func{grad}f)%
\end{array}%
\right\}  \nonumber \\
&=&\sum\limits_{i=1}^{m}\left\{
\begin{array}{c}
e_{i}(e_{i}(f))\tau (\varphi )+2e_{i}(f)\nabla _{e_{i}}^{\varphi }\,\tau
(\varphi )+f\nabla _{e_{i}}^{\varphi }\nabla _{e_{i}}^{\varphi }\,\tau
(\varphi )+\nabla _{e_{i}}^{\varphi }\nabla _{e_{i}}^{\varphi }d\varphi
(\func{grad}f) \\
-\left( \nabla _{e_{i}}^{M}e_{i}\right) (f)\tau (\varphi )-f\nabla _{\nabla
_{e_{i}}^{M}e_{i}}^{\varphi }\tau (\varphi )-\nabla _{\nabla
_{e_{i}}^{M}e_{i}}^{\varphi }d\varphi (\func{grad}f)%
\end{array}%
\right\}  \nonumber \\
&=&\left( \Delta f\right) \tau (\varphi )+2\nabla _{\func{grad}f}^{N}\,\tau
(\varphi )-f\Delta ^{\varphi }\left( \tau (\varphi )\right) -\Delta
^{\varphi }\left( \func{grad}f\right) .  \label{hp3}
\end{eqnarray}%

Since the tension field of $\varphi $ is given by $\tau (\varphi )=mH\xi $,
we have%
\begin{equation}
\left( \Delta f\right) \tau (\varphi )=mH\left( \Delta f\right) \xi ,
\label{hp4}
\end{equation}%
\begin{eqnarray}
\nabla _{\func{grad}f}^{N}\,\tau (\varphi ) &=&\nabla _{\func{grad}%
f}^{N}\,(mH\xi )=m\func{grad}f(H)\xi -mHA(\func{grad}f)  \nonumber \\
&=&m\left\langle \func{grad}f,\func{grad}H\right\rangle \xi -mHA(\func{grad}%
f),  \label{hp5}
\end{eqnarray}%
\begin{eqnarray}
\Delta ^{\varphi }\left( \tau (\varphi )\right)
&=&-\sum\limits_{i=1}^{m}\left\{ \nabla _{e_{i}}^{\varphi }\nabla
_{e_{i}}^{\varphi }(mH\xi )-\nabla _{\nabla _{e_{i}}^{M}e_{i}}^{\varphi
}(mH\xi )\right\}  \nonumber \\
&=&-m\sum\limits_{i=1}^{m}\left\{ \nabla _{e_{i}}^{\varphi }\left(
e_{i}(H)\xi +H\nabla _{e_{i}}^{N}\xi \right) -\left( \nabla
_{e_{i}}^{M}e_{i}\right) \left( H\right) \xi -H\nabla _{\nabla
_{e_{i}}^{M}e_{i}}^{N}\xi \right\}  \nonumber \\
&=&-m\sum\limits_{i=1}^{m}\left\{ e_{i}(e_{i}(H))\xi +2e_{i}(H)\nabla
_{e_{i}}^{N}\xi +H\nabla _{e_{i}}^{N}\nabla _{e_{i}}^{N}\xi -\left( \nabla
_{e_{i}}^{M}e_{i}\right) \left( H\right) \xi -H\nabla _{\nabla
_{e_{i}}^{M}e_{i}}^{N}\xi \right\}  \nonumber \\
&=&-m\left( \Delta H\right) \xi +2mA(\func{grad}H)+mH\Delta ^{\varphi
}\left( \xi \right) .  \label{hp6}
\end{eqnarray}%

By using (\ref{hp4}), (\ref{hp5}) and (\ref{hp6}) in (\ref{hp3}) we obtain%
\begin{eqnarray}
\func{trace}(\nabla ^{\varphi })^{2}\tau _{f\,}\left( \varphi \right) &=&mH\left(
\Delta f\right) \xi +2m\left\langle \func{grad}f,\func{grad}H\right\rangle
\xi -2mHA(\func{grad}f) \nonumber \\
&&+mf\left( \Delta H\right) \xi -2mfA(\func{grad}H)-mfH\Delta
^{\varphi }\left( \xi \right) -\Delta ^{\varphi }\left( \func{grad}f\right).   \label{hp7}
\end{eqnarray}%

As a second step, we compute:%
\begin{eqnarray}
\func{trace} R^{N}\left( \tau _{f\,}\left( \varphi \right) ,d\varphi \right)
d\varphi &=&\sum\limits_{i=1}^{m}R^{N}(f\,\tau (\varphi )+d\varphi
(\func{grad}f),d\varphi (e_{i}))d\varphi (e_{i}) \nonumber  \\
&=&f\,\sum\limits_{i=1}^{m}R^{N}(\tau (\varphi
),e_{i})e_{i}+\sum\limits_{i=1}^{m}R^{N}(\func{grad}f,e_{i})e_{i},  \label{hp8}
\end{eqnarray}%
which implies
\begin{equation}
\func{trace}R^{N}\left( \tau _{f\,}\left( \varphi \right) ,d\varphi \right)
d\varphi =mfH\,\sum\limits_{i=1}^{m}R^{N}(\xi
,e_{i})e_{i}+\sum\limits_{i=1}^{m}R^{N}(\func{grad}f,e_{i})e_{i}.  \label{hp9}
\end{equation}%

Also
\begin{eqnarray}
\nabla _{\func{grad}f}^{\varphi }\tau _{f}\,\left( \varphi \right) &=&\nabla
_{\func{grad}f}^{\varphi }\left( f\,\tau (\varphi )+d\varphi (\func{grad}f)\right)
\nonumber \\
&=&\nabla _{\func{grad}f}^{N}\left( f\,\tau (\varphi )\right) +\nabla _{%
\func{grad}f}^{N}\func{grad}f  \nonumber \\
&=&\func{grad}f(f)\tau (\varphi )+f\nabla _{\func{grad}f}^{N}\,\tau (\varphi
)+\nabla _{\func{grad}f}^{N}\func{grad}f,  \label{hp10}
\end{eqnarray}%
which gives%
\begin{eqnarray}
\nabla _{\func{grad}f}^{\varphi }\tau _{f}\,\left( \varphi \right)
&=&mH\left\langle \func{grad}f,\func{grad}f\right\rangle \xi +mf\left\langle
\func{grad}f,\func{grad}H\right\rangle \xi  \nonumber \\
&&-mfHA(\func{grad}f)+\frac{1}{2}\func{grad}(\left\vert \func{grad}%
f\right\vert ^{2})+B(\func{grad}f,\func{grad}f).  \label{hp11}
\end{eqnarray}%

By using (\ref{hp7}), (\ref{hp9}) and (\ref{hp11}) in (\ref{hp2}) we obtain
the bi-$f$-tension field of $\varphi $:%
\begin{eqnarray}
\tau _{f,2}(\varphi ) &=&-mfH\left( \Delta f\right) \xi -3mf\left\langle
\func{grad}f,\func{grad}H\right\rangle \xi +3mfHA(\func{grad}f)  \nonumber
\\
&&-mf^{2}\left( \Delta H\right) \xi +2mf^{2}A(\func{grad}H)+mf^{2}H
\Delta ^{\varphi }\left(\xi \right)  \nonumber \\
&&+f\Delta ^{\varphi }(\func{grad}f)-mH\left\langle \func{grad}f,\func{grad}%
f\right\rangle \xi  \nonumber \\
&&-\frac{1}{2}\func{grad}(\left\vert \func{grad}f\right\vert ^{2})-B(\func{grad}f,%
\func{grad}f)  \nonumber \\
&&-mf^{2}H\sum\limits_{i=1}^{m}R^{N}(\xi
,e_{i})e_{i}-f\sum\limits_{i=1}^{m}R^{N}(\func{grad}f,e_{i})e_{i}. \label{hp12}
\end{eqnarray}



The tangential component of $\Delta ^{\varphi }(\func{grad}f)$ can be calculated by%
\begin{eqnarray}
\left( \Delta ^{\varphi }(\func{grad}f)\right) ^{\top }
&=&-\sum_{i,k=1}^{m}\left\langle \nabla _{e_{i}}^{N}\nabla _{e_{i}}^{N}\func{%
grad}f-\nabla _{\nabla _{e_{i}}^{M}e_{i}}^{N}\func{grad}f,e_{k}\right\rangle
e_{k}  \nonumber \\
&=&-\sum_{i,k=1}^{m}\left\langle \nabla _{e_{i}}^{N}\left( \nabla
_{e_{i}}^{M}\func{grad}f+B(\func{grad}f,e_{i}\right) ),e_{k}\right\rangle
e_{k}  \nonumber \\
&&-\sum_{i,k=1}^{m}\left\langle -\nabla _{\nabla _{e_{i}}^{M}e_{i}}^{M}\func{%
grad}f-B(\nabla _{e_{i}}^{M}e_{i},\func{grad}f),e_{k}\right\rangle e_{k}
\nonumber \\
&=&-\sum_{i,k=1}^{m}\left\langle \nabla _{e_{i}}^{M}\nabla _{e_{i}}^{M}\func{%
grad}f+B(\nabla _{e_{i}}^{M}\func{grad}f,e_{i})+\nabla _{e_{i}}^{N}B(\func{%
grad}f,e_{i}),e_{k}\right\rangle e_{k} \nonumber  \\
&&-\sum_{i,k=1}^{m}\left\langle -\nabla _{\nabla _{e_{i}}^{M}e_{i}}^{M}\func{%
grad}f-B(\nabla _{e_{i}}^{M}e_{i},\func{grad}f),e_{k}\right\rangle e_{k}
\nonumber \\
&=&-\sum_{i,k=1}^{m}\left\langle \nabla _{e_{i}}^{M}\nabla _{e_{i}}^{M}\func{%
grad}f+b(\nabla _{e_{i}}^{M}\func{grad}f,e_{i})\xi+\nabla _{e_{i}}^{N}\left( b(%
\func{grad}f,e_{i})\xi \right) ,e_{k}\right\rangle e_{k}  \nonumber \\
&&-\sum_{i,k=1}^{m}\left\langle -\nabla _{\nabla _{e_{i}}^{M}e_{i}}^{M}\func{%
grad}f-b(\nabla _{e_{i}}^{M}e_{i},\func{grad}f)\xi,e_{k}\right\rangle e_{k}
\nonumber \\
&=&-\sum_{i,k=1}^{m}\left\langle \nabla _{e_{i}}^{M}\nabla _{e_{i}}^{M}\func{%
grad}f-\nabla _{\nabla _{e_{i}}^{M}e_{i}}^{M}\func{grad}f,e_{k}\right\rangle
e_{k}  \nonumber \\
&&-\sum_{i,k=1}^{m}\left\langle b(\nabla _{e_{i}}^{M}\func{grad}f,e_{i})\xi
+\nabla _{e_{i}}^{N}\left( b(\func{grad}f,e_{i})\xi \right) -b(\nabla
_{e_{i}}^{M}e_{i},\func{grad}f)\xi ,e_{k}\right\rangle e_{k}  \nonumber \\
&=&-\sum_{i,k=1}^{m}\left\langle \nabla _{e_{i}}^{M}\nabla _{e_{i}}^{M}\func{%
grad}f-\nabla _{\nabla _{e_{i}}^{M}e_{i}}^{M}\func{grad}f+b(\nabla
_{e_{i}}^{M}\func{grad}f,e_{i})\xi ,e_{k}\right\rangle e_{k}  \nonumber \\
&&-\sum_{i,k=1}^{m}\left\langle e_{i}\left( b(\func{grad}f,e_{i})\right) \xi
-b(\func{grad}f,e_{i})Ae_{i}-b(\nabla _{e_{i}}^{M}e_{i},\func{grad}f)\xi
,e_{k}\right\rangle e_{k}  \nonumber \\
&=&\Delta (\func{grad}f)+\sum_{i,k=1}^{m}b(\func{grad}f,e_{i})\left\langle
Ae_{i},e_{k}\right\rangle e_{k},  \label{hp12i}
\end{eqnarray}%
which implies%
\begin{equation}
\left( \Delta ^{\varphi }(\func{grad}f)\right) ^{\top }=\Delta (\func{grad}%
f)+A^{2}(\func{grad}f).  \label{hp12ii}
\end{equation}%

The normal component of $\Delta ^{\varphi }(\func{grad}f)$ can be calculated by%
\begin{eqnarray}
\left( \Delta ^{\varphi }(\func{grad}f)\right) ^{\bot }
&=&-\sum_{i=1}^{m}\left\langle \nabla _{e_{i}}^{N}\nabla _{e_{i}}^{N}\func{%
grad}f-\nabla _{\nabla _{e_{i}}^{M}e_{i}}^{N}\func{grad}f,\xi \right\rangle
\xi  \nonumber \\
&=&-\sum_{i=1}^{m}\left\langle
\begin{array}{c}
\nabla _{e_{i}}^{N}\left( \nabla _{e_{i}}^{M}\func{grad}f+B(\func{grad}%
f,e_{i}\right) ) \\
-\nabla _{\nabla _{e_{i}}^{M}e_{i}}^{M}\func{grad}f-B(\nabla
_{e_{i}}^{M}e_{i},\func{grad}f)%
\end{array}%
,\xi \right\rangle \xi  \nonumber \\
&=&-\sum_{i=1}^{m}\left\langle
\begin{array}{c}
\nabla _{e_{i}}^{M}\nabla _{e_{i}}^{M}\func{grad}f+B(\nabla _{e_{i}}^{M}%
\func{grad}f,e_{i}) \\
+\nabla _{e_{i}}^{N}B(\func{grad}f,e_{i})-\nabla _{\nabla
_{e_{i}}^{M}e_{i}}^{M}\func{grad}f-B(\nabla _{e_{i}}^{M}e_{i},\func{grad}f)%
\end{array}%
,\xi \right\rangle \xi  \nonumber \\
&=&-\sum_{i=1}^{m}\left\langle
\begin{array}{c}
\nabla _{e_{i}}^{M}\nabla _{e_{i}}^{M}\func{grad}f+b(\nabla _{e_{i}}^{M}%
\func{grad}f,e_{i})\xi \\
+\nabla _{e_{i}}^{N}\left( b(\func{grad}f,e_{i})\xi \right) -\nabla _{\nabla
_{e_{i}}^{M}e_{i}}^{M}\func{grad}f-b(\nabla _{e_{i}}^{M}e_{i},\func{grad}f)\xi%
\end{array}%
,\xi \right\rangle \xi \nonumber  \\
&=&-\sum_{i=1}^{m}\left\langle
\begin{array}{c}
\nabla _{e_{i}}^{M}\nabla _{e_{i}}^{M}\func{grad}f-\nabla _{\nabla
_{e_{i}}^{M}e_{i}}^{M}\func{grad}f+b(\nabla _{e_{i}}^{M}\func{grad}%
f,e_{i})\xi \\
+\nabla _{e_{i}}^{N}\left( b(\func{grad}f,e_{i})\xi \right) -b(\nabla
_{e_{i}}^{M}e_{i},\func{grad}f)\xi%
\end{array}%
,\xi \right\rangle \xi  \nonumber \\
&=&-\sum_{i=1}^{m}\left\langle
\begin{array}{c}
\nabla _{e_{i}}^{M}\nabla _{e_{i}}^{M}\func{grad}f-\nabla _{\nabla
_{e_{i}}^{M}e_{i}}^{M}\func{grad}f+b(\nabla _{e_{i}}^{M}\func{grad}%
f,e_{i})\xi \\
+e_{i}\left( b(\func{grad}f,e_{i})\right) \xi -b(\func{grad}%
f,e_{i})Ae_{i}-b(\nabla _{e_{i}}^{M}e_{i},\func{grad}f)\xi%
\end{array}%
,\xi \right\rangle \xi  \nonumber \\
&=&-\sum_{i=1}^{m}\left\langle b(\nabla _{e_{i}}^{M}\func{grad}f,e_{i})\xi
+e_{i}\left( b(\func{grad}f,e_{i})\right) \xi -b(\nabla _{e_{i}}^{M}e_{i},%
\func{grad}f)\xi ,\xi \right\rangle \xi  \nonumber \\
&=& -\sum_{i=1}^{m}\{b(\nabla _{e_{i}}^{M}\func{grad}f,e_{i})
+e_{i}\left( b(\func{grad}f,e_{i})\right) -b(\nabla _{e_{i}}^{M}e_{i},%
\func{grad}f)\} \xi. \label{hp12iii}
\end{eqnarray}


The tangential component of $\Delta ^{\varphi }(\xi)$ can be calculated by%
\begin{eqnarray}
\left( \Delta ^{\varphi }(\xi) \right) ^{\top }
&=&-\sum_{i,k=1}^{m}\left\langle \nabla _{e_{i}}^{N}\nabla _{e_{i}}^{N}\xi
-\nabla _{\nabla _{e_{i}}^{M}e_{i}}^{N}\xi ,e_{k}\right\rangle e_{k}
\nonumber \\
&=&\sum_{i,k=1}^{m}\left\langle \nabla _{e_{i}}^{N}Ae_{i}-A(\nabla
_{e_{i}}^{M}e_{i}),e_{k}\right\rangle e_{k}  \nonumber \\
&=&\sum_{i,k=1}^{m}\{e_{i}\left\langle Ae_{i},e_{k}\right\rangle
-\left\langle Ae_{i},\nabla _{e_{i}}^{M}e_{k}\right\rangle -\left\langle
A(\nabla _{e_{i}}^{M}e_{i}),e_{k}\right\rangle \}e_{k} \nonumber  \\
&=&\sum_{i,k=1}^{m}\{e_{i}b(e_{i},e_{k})-b(e_{i},\nabla
_{e_{i}}^Me_{k})-b(\nabla _{e_{i}}^Me_{i},e_{k})\}e_{k}  \nonumber \\
&=&\sum_{i,k=1}^{m}\{(\nabla _{e_{i}}^Mb)(e_{k},e_{i})\}e_{k}.  \label{hp12a}
\end{eqnarray}%

By Codazzi-Mainardi equation, we have%
\begin{equation}
\sum_{i=1}^{m}((\nabla _{e_{i}}^Mb)(e_{k},e_{i})-(\nabla
_{e_{k}}^Mb)(e_{i},e_{i}))=-\sum_{i=1}^{m}\left\langle
R^{N}(e_{i},e_{k})e_{i},\xi \right\rangle =Ric^{N}(\xi ,e_{k}).
\label{hp12b}
\end{equation}%

Putting the last equation into (\ref{hp12a}) we get%
\begin{eqnarray}
\left( \Delta ^{\varphi }(\xi )\right) ^{\top } &=&\sum_{i,k=1}^{m}\{(\nabla
_{e_{i}}^Mb)(e_{k},e_{i})\}e_{k} \nonumber \\
&=&\sum_{k=1}^{m}\{\sum_{i=1}^{m}(\nabla _{e_{k}}^Mb)(e_{i},e_{i})+Ric^{N}(\xi ,e_{k})\}e_{k}
\nonumber \\
&=&m\,{\func{grad}}H+\sum_{k=1}^{m}Ric^{N}(\xi ,e_{k})e_{k}. \label{p}
\end{eqnarray}%

The normal component of $\Delta ^{\varphi }(\xi)$ can be calculated by
\begin{eqnarray}
\left( \Delta ^{\varphi }(\xi) \right) ^{\bot }
&=&-\sum_{i=1}^{m}\left\langle \nabla _{e_{i}}^{N}\nabla _{e_{i}}^{N}\xi
-\nabla _{\nabla _{e_{i}}^{M}e_{i}}^{N}\xi ,\xi \right\rangle \xi  \nonumber \\
&=&-\sum_{i=1}^{m}\left\langle \nabla _{e_{i}}^{N}\nabla _{e_{i}}^{N}\xi
,\xi \right\rangle \xi \nonumber \\
&=&\sum_{i=1}^{m}\left\langle \nabla _{e_{i}}^{N}\xi ,\nabla _{e_{i}}^{N}\xi
\right\rangle \xi.  \label{hp13}
\end{eqnarray}%

On the other hand
\begin{eqnarray}
\sum_{i=1}^{m}\left\langle \nabla _{e_{i}}^{N}\xi ,\nabla _{e_{i}}^{N}\xi
\right\rangle &=&\sum_{i,j=1}^{m}\left\langle \nabla _{e_{i}}^{N}\xi
,\left\langle \nabla _{e_{i}}^{N}\xi ,e_{j}\right\rangle e_{j}\right\rangle
\nonumber \\
&=&\sum_{i,j=1}^{m}\left\langle \nabla _{e_{i}}^{N}\xi ,e_{j}\right\rangle
^{2}  \nonumber \\
&=&\sum_{i,j=1}^{m}\left\langle Ae_{i},e_{j}\right\rangle ^{2}  \nonumber \\
&=&\left\vert A\right\vert ^{2},  \label{hp14}
\end{eqnarray}%
which implies together with (\ref{hp13})%
\begin{equation}
\left( \Delta ^{\varphi }(\xi) \right) ^{\bot }=\left\vert A\right\vert ^{2}\xi.  \label{hp15}
\end{equation}


The tangential and the normal components of the curvature terms are
\begin{equation}
\sum\limits_{i,k=1}^{m}\left\langle R^{N}(\xi
,e_{i})e_{i},e_{k}\right\rangle e_{k}=\sum_{k=1}^{m}Ric^{N}(\xi ,e_{k})e_{k}=(Ric^{N}(\xi))^{\top },  \label{hp16}
\end{equation}%
\begin{equation}
\sum\limits_{i=1}^{m}\left\langle R^{N}(\xi ,e_{i})e_{i},\xi \right\rangle \xi
=Ric^{N}(\xi ,\xi )\xi,  \label{hp17}
\end{equation}%
\begin{equation}
\sum\limits_{i,k=1}^{m}\left\langle R^{N}(\func{grad}f,e_{i})e_{i},e_{k}\right\rangle
e_{k}=\sum_{k=1}^{m}Ric^{N}(\func{grad}f,e_{k})e_{k}=(Ric^{N}(\func{grad}f))^{\top },  \label{hp18}
\end{equation}%
\begin{equation}
\sum\limits_{i=1}^{m}\left\langle R^{N}(\func{grad}f,e_{i})e_{i},\xi \right\rangle \xi
=Ric^{N}(\func{grad}f,\xi )\xi.  \label{hp19}
\end{equation}%

By collecting all the tangential and normal components of the bi-$f$%
-tension field separately, we have%
\begin{eqnarray}
\left[ \tau _{f,2}(\varphi )\right] ^{\top } &=&3mfHA(\func{grad}f)+2mf^{2}A(%
\func{grad}H)  \nonumber \\
&&+m^{2}f^{2}H{\func{grad}}H+f\Delta (\func{grad}f)  \nonumber \\
&&+fA^{2}(\func{grad}f)-\frac{1}{2}\func{grad}(\left\vert \func{grad}%
f\right\vert ^{2}) \nonumber  \\
&&-mf^2H(Ric^{N}(\xi))^{\top }-f\left( Ric^{N}(\func{grad}%
f)\right) ^{\top }  \label{hp20}
\end{eqnarray}%
and
\begin{eqnarray}
\left[ \tau _{f,2}(\varphi )\right] ^{\bot } &=&\{-mfH\left( \Delta f\right)
-3mf\left\langle \func{grad}f,\func{grad}H\right\rangle  \nonumber \\
&&-mf^{2}\left( \Delta H\right) +mf^{2}H\left\vert A\right\vert ^{2}+f
(\Delta ^{\varphi }\left(\func{grad}f\right)) ^{\bot }  \nonumber \\
&&-mH\left\vert \func{grad}f\right\vert ^{2}-b(\func{grad}f,\func{grad}f)
\nonumber \\
&&-mf^{2}HRic^{N}(\xi ,\xi )-fRic^{N}(\func{grad}f,\xi )\}\xi.  \label{hp21}
\end{eqnarray}%

Then we have

\begin{theorem}
Let $(N,h)$ be an $(m+1)$-dimensional Riemannian manifold and $\varphi :M\rightarrow N$ be an isometric immersion of
codimension-one with mean curvature vector $\eta $ $=H\xi $. Then $\varphi $
is a bi-$f$-harmonic map if and only if%
\begin{eqnarray}
0 &=&3mfHA(\func{grad}f)+2mf^{2}A(\func{grad}H)+m^{2}f^{2}H{\func{grad}}%
H  \nonumber \\
&&+f\Delta (\func{grad}f)+fA^{2}(\func{grad}f)-\frac{1}{2}\func{grad}%
(\left\vert \func{grad}f\right\vert ^{2}) \nonumber \\
&&-mf^2H(Ric^{N}(\xi))^{\top }-f\left( Ric^{N}(\func{grad}%
f)\right) ^{\top }  \label{hp22}
\end{eqnarray}%
and
\begin{eqnarray}
0 &=&-mfH\left( \Delta f\right) -3mf\left\langle \func{grad}f,\func{grad}%
H\right\rangle -mf^{2}\left( \Delta H\right) +mf^{2}H\left\vert A\right\vert
^{2}+f\left( \Delta ^{\varphi }(\func{grad}f)\right) ^{\bot }  \nonumber \\
&&-mH\left\vert \func{grad}f\right\vert ^{2}-b(\func{grad}f,\func{grad}%
f)-mf^{2}HRic^{N}(\xi ,\xi )-fRic^{N}(\func{grad}f,\xi ),  \label{hp23}
\end{eqnarray}%
where $Ric^{N}$ denotes also the Ricci operator of the ambient space, $A$ is the shape operator of the hypersurface with respect to
the unit normal vector field $\xi $, $\Delta $ and $\func{grad}$ are the Laplace and
the gradient operator of the hypersurface, respectively, and $\Delta ^{\varphi }$ is the rough Laplace operator on sections of $\varphi ^{-1}TN$.
\end{theorem}

\begin{theorem}
Let $M$ be a constant mean curvature hypersurface in an $(m+1)$-dimensional Riemannian manifold $%
N.$ Then $M$ is a bi-$f$-harmonic submanifold if and only if
\begin{eqnarray}
mf^2H(Ric^{N}(\xi))^{\top }+f\left( Ric^{N}(\func{grad}%
f)\right) ^{\top } &=&3mfHA(\func{grad}f)+f\Delta (\func{grad}f)
\nonumber \\
&&+fA^{2}(\func{grad}f)-\frac{1}{2}\func{grad}(\left\vert \func{grad}%
f\right\vert ^{2})  \label{hp24}
\end{eqnarray}%
and%
\begin{eqnarray}
mf^{2}HRic^{N}(\xi ,\xi )+fRic^{N}(\func{grad}f,\xi ) &=&-mfH\left( \Delta f\right)
+mf^{2}H\left\vert A\right\vert ^{2}+f\left( \Delta ^{\varphi }(\func{grad}%
f)\right) ^{\bot }  \nonumber \\
&&-mH\left\vert \func{grad}f\right\vert ^{2}-b(\func{grad}f,\func{grad}f). \label{hp25}
\end{eqnarray}
\end{theorem}

Then we have

\begin{corollary}
Let $M$ be a constant mean curvature hypersurface in an $(m+1)$-dimensional Ricci flat Riemannian
manifold $N.$ Then $M$ is a bi-$f$-harmonic submanifold if and only if%
\begin{equation}
fA^{2}(\func{%
grad}f)+3mfHA(\func{grad}f)+f\Delta (\func{grad}f)-\frac{1}{2}\func{grad}(\left\vert \func{grad}f\right\vert ^{2})=0  \label{hp24a}
\end{equation}%
and%
\begin{equation}
mfH\left( \Delta f\right) +mH\left\vert
\func{grad}f\right\vert ^{2}-mf^{2}H\left\vert A\right\vert ^{2}+b(\func{grad}f,\func{grad}f)-f\left( \Delta
^{\varphi }(\func{grad}f)\right) ^{\bot }=0.  \label{hp24b}
\end{equation}
\end{corollary}

\begin{corollary}
Let $M$ be a hypersurface in an $(m+1)$-dimensional Einstein space $N$. Then $M$ is a bi-$%
f$-harmonic submanifold if and only if
\begin{eqnarray}
f\frac{r}{m+1}\func{grad}f &=&3mfHA(\func{grad}f)+f\Delta (\func{grad}f)
\nonumber \\
&&+fA^{2}(\func{grad}f)-\frac{1}{2}\func{grad}(\left\vert \func{grad}%
f\right\vert ^{2})  \label{hp26}
\end{eqnarray}%
and%
\begin{eqnarray}
mf^{2}H\frac{r}{m+1} &=&-mfH\left( \Delta f\right) +mf^{2}H\left\vert
A\right\vert ^{2}+f\left( \Delta ^{\varphi }(\func{grad}f)\right) ^{\bot }
\nonumber \\
&&-mH\left\vert \func{grad}f\right\vert ^{2}-b(\func{grad}f,\func{grad}f),
\label{hp27}
\end{eqnarray}%
where $r$ is the scalar curvature of the ambient space.
\end{corollary}

Since an $(m+1)$-dimensional space of constant sectional curvature $c$
is an Einstein space with scalar curvature $r=m(m+1)c,$ by using (\ref%
{hp26}) and (\ref{hp27}) we have

\begin{corollary}
Let $M$ be a hypersurface in an $(m+1)$-dimensional space $N$ of constant sectional curvature $c$. Then $M$ is a bi-$%
f$-harmonic submanifold if and only if
\begin{eqnarray}
mcf\func{grad}f&=&3mfHA(\func{grad}f) +f\Delta (\func{grad}f)  \nonumber
\\
&&+fA^{2}(\func{grad}f)-\frac{1}{2}\func{grad}(\left\vert \func{grad}%
f\right\vert ^{2}) \label{hp28}
\end{eqnarray}%
and%
\begin{eqnarray}
m^{2}cf^{2}H &=&-mfH\left( \Delta f\right) +mf^{2}H\left\vert A\right\vert
^{2}+f\left( \Delta ^{\varphi }(\func{grad}f)\right) ^{\bot } \nonumber  \\
&&-mH\left\vert \func{grad}f\right\vert ^{2}-b(\func{grad}f,\func{grad}f).
\label{hp29}
\end{eqnarray}
\end{corollary}


\bigskip

\small{
{\bf Selcen Y\"{u}ksel Perkta\c{s}} \\
Department of Mathematics, Faculty of Arts and Sciences \\
Ad\i yaman University,
02040, Ad\i yaman, Turkey \\
e-mail: sperktas@adiyaman.edu.tr

\smallskip

{\bf Adara Monica Blaga} \\
Faculty of Mathematics and Computer Science, Department of Mathematics \\
West University of Timi\c{s}oara,
300223, Timi\c{s}oara, Rom\^{a}nia \\
e-mail: adarablaga@yahoo.com

\smallskip

{\bf Feyza Esra Erdo\u{g}an} \\
Faculty of Education, Department of Elementary Education \\
Ad\i yaman University,
02040, Ad\i yaman, Turkey \\
e-mail: ferdogan@adiyaman.edu.tr

\smallskip

{\bf Bilal Eftal Acet} \\
Department of Mathematics \\
Ad\i yaman University,
02040, Ad\i yaman, Turkey \\
e-mail: eacet@adiyaman.edu.tr
}

\end{document}